\documentclass[11pt, 14paper,reqno]{amsart}
\usepackage{amsmath,amsfonts,amssymb}
\usepackage{bigints}
\usepackage{longtable}
\theoremstyle{plain}
\newtheorem{theorem}{Theorem}

\theoremstyle{proof}
\theoremstyle{definition}

\newtheorem{prop}{Proposition}

\theoremstyle{remark}
\newtheorem{remark}{Remark}
\theoremstyle{lamma}

\renewcommand{\Re}{{\mathrm Re \,}}

\newcommand{\F}{{\mathcal{F}}}
\newcommand{\R}{{\mathfrak{R} }}
\newcommand{\z}{{\zeta}}
\newcommand{\ka}{{\kappa}}

\numberwithin{equation}{section}
\numberwithin{lemma}{section}
\numberwithin{theorem}{section}
\numberwithin{remark}{section}
\numberwithin{prop}{section}

\usepackage{amsmath}
\usepackage{amsfonts}   
\usepackage{amssymb}
\usepackage{amssymb, amsmath, amsthm}
\usepackage[breaklinks]{hyperref}
\usepackage{cases}
\theoremstyle{thmrm} 

\newtheorem{conjecture}{Conjecture}
\numberwithin{conjecture}{section}

\newtheorem{thma}{Theorem}

\newenvironment{dedication}
    {\vspace{3ex}\begin{quotation}\begin{center}\begin{em}}
    {\par\end{em}\end{center}\end{quotation}}

\begin{document}
\title[Values of Dedekind zeta functions]{An analogue of Wilton's formula and values of Dedekind zeta functions}
\author{Soumyarup Banerjee, Kalyan Chakraborty and Azizul Hoque}
\address{Soumyarup Banerjee @Department of Mathematics, IIT Gandhinagar, 
Palaj, Gandhinagar-382355, Gujarat, India.}
\email{soumya.tatan@gmail.com}

\address{Kalyan Chakraborty @Kerala School of Mathematics, 
Kunnamangalam-673571, Kozhikode, Kerala, India.}
\email{ kalychak@ksom.res.in}
\address{Azizul Hoque @Harish-Chandra Research Institute,
Chhatnag Road, Jhunsi,  Allahabad 211 019, India.}
\email{ahoque.ms@gmail.com}
\maketitle
\begin{dedication}
To Prof. Bruce C. Berndt on his $80^{th}$ birthday with great respect
\end{dedication}
\begin{abstract}
J. R. Wilton obtained an expression for the product of two Riemann zeta functions. This expression played a crucial role to find the approximate functional equation for the product of two Riemann zeta functions in the critical region. We find analogous expressions for the product of two Dedekind zeta functions and then use these expressions to find some expressions for Dedekind zeta values attached to arbitrary real as well as quadratic number fields at any positive integer. 
\end{abstract}
\tableofcontents
\section{Introduction}
G. H. Hardy and J. E. Littlewood \cite{HL1921} obtained an approximate functional equation for Riemann zeta function which is particularly useful in studying the behaviour of zeta function in the critical strip; especially on the critical line. Later, in the same way of thought, J. R. Wilton \cite{WI30} in 1930 gave an approximate functional equation for the product of two Riemann zeta functions in the critical region. The main ingredient used to arrive on to the above mentioned approximate functional equation is the following formula:   
\begin{thma}[Wilton's formula]\label{Wilton}
Let $\R (u)>-1, ~ \R (v)>-1,~ \R (u+v)>0$ and $u+v\ne 2$, then
\begin{equation*}
\begin{aligned}
\z(u)&\z(v)= \z(u+v-1)\bigg(\frac{1}{u-1}+\frac{1}{v-1}\bigg)+2(2\pi)^{u-1}\sum_{n=1}^\infty\sigma_{1-u-v}(n)n^{u-1}u\\
\times& \int_{2n\pi}^\infty x^{-u-1}\sin x \  \rm{d}x  + 2(2\pi)^{v-1}\sum_{n=1}^\infty\sigma_{1-u-v}(n)n^{v-1}v \int_{2n\pi}^\infty x^{-v-1}\sin x \ \rm{d}x,
\end{aligned}
\end{equation*}
where $\sigma_z(n)$, $z$ being an arbitrary (complex) number, denotes the sum of $d^z$ over all the divisors $d$ of $n$.
\end{thma}
Wilton's approach to prove Theorem \ref{Wilton} involves techniques of real and complex analysis. Recently, D. Banerjee and J. Mehta \cite{BM14} gave an alternative proof of Theorem \ref{Wilton} using number theoretic approach, and a similar result for Hurwitz zeta functions appeared in \cite{WB2017}. In this paper, we obtain a result analogous to that of Theorem \ref{Wilton} in number fields. More precisely, we find expressions for the product of two Dedekind zeta functions attached to arbitrary real (resp. imaginary) quadratic number fields. These expressions involve Meijor G-function as its error term. We then use these expressions to find some formulae for Dedekind zeta values attached to quadratic number fields at any positive integer. The ultimate goal is to find certain formulae for calculating special values of Dedekind zeta functions, which we are unable to find at this point. To the best of our knowledge these expressions are new and are of some interests.
\subsection{Preliminaries}\hspace{10cm}$ $
We shall begin by recalling some definitions and fix up notations to be used in this paper. Let $K$ be a number field of degree $d$ with signature $(r_1, r_2)$ (i. e., $d=r_1+2r_2$) and $\mathcal{O}_K$ it's ring of integers.
Let $N(\mathfrak{a})$ be the norm of the ideal $\mathfrak{a}\subseteq \mathcal{O}_K$. Then the Dedekind zeta function attached to $K$ is defined by 
$$
\z_K(s)=\sum_{\mathfrak{a}\subset\mathcal{O}_K}\frac{1}{N(\mathfrak{a})^s}=\prod_{\mathfrak{p}\subset \mathcal{O}_K}\bigg(1-\frac{1}{N(\mathfrak{p})^s}\bigg)^{-1},\hspace*{4mm} \R (s)>1, 
$$
where $\mathfrak{a}$ and $\mathfrak{p}$ run over the non-zero ideals and prime ideals of $\mathcal{O}_K$ respectively.
If $v_K(m)$ denotes the number of non-zero integral ideals in $\mathcal{O}_K$ with norm $m$, then $\z_K$ can also be expressed as 
$$
\z_K(s)=\sum_{m=1}^\infty \frac{v_K(m)}{m^s}.$$
Let 
$$
\sigma'_{m, v_K}(n):=\sum_{t\mid n}t^mv_K(t)v_K(\frac{n}{t}),
$$
for any integer $m$, and any positive integer $n$.
\\
The Meijer G-function is a sum of hypergeometric 
functions each of which is usually an entire function. It is defined by the following line integral in the complex plane (see \cite{N62})
\begin{equation}\label{eq1}
\begin{aligned}
G^{m, \ n}_{p, \ q}\bigg(\begin{matrix}
a_1, \cdots, a_p \\
b_1, \cdots, b_q
\end{matrix} \ \bigg|\ z\bigg)=\frac{1}{2\pi i}\underset{{(C)}}{\bigints} \frac{\prod\limits_{j=1}^m\Gamma(b_j-s)\prod\limits_{j=1}^n\Gamma(1-a_j+s)}{\prod\limits_{j=m+1}^q\Gamma(1-b_j+s)\prod\limits_{j=n+1}^p\Gamma(a_j-s)}z^s \rm{d}s .
\end{aligned}\nonumber
\end{equation}
The poles of the integrand must be simple and those of $\Gamma(b_j-s), j=1, \cdots, m$, must lie on one side of the contour $C$ and those of
$\Gamma(1-a_j+s), j=1, \cdots, n$, must lie on the other side of $C$. The integral then converges for $|\arg z| < \delta \pi$,  where
$$\delta = m+n - \frac{1}{2}(p+q).$$
The integral additionally converges for $|\arg z|= \delta \pi$ if $(q-p)(\Re(s) + 1/2) > \Re (v) + 1$, where 
$$v = \sum_{j=1}^{q}b_j - \sum_{j=1}^{p}a_j.$$ 

\subsection{Statement of results}\hspace{10cm}$ $
Let $D_K,\ h_K, \ R_K$ and $w_{K}$  be denote the discriminant, the class number, the regulator and  the number of roots of unity in the number field $K$, respectively.  We further define the following auxiliary functions: 
\begin{eqnarray*}
Z(z):=&&4(2\pi)^{2(z-1)}D_K ^{\frac{1-2z}{2}}\sum_{m=1}^\infty \sigma'_{1-u-v, v_K}(m)m^{z-1}\bigg\{\cos \pi z~~ G^{0, 3}_{3, 1}\bigg(\begin{matrix}
1, z,  z \\
0, -
\end{matrix}\ \bigg|\ \frac{-D_K}{4\pi ^2m}\bigg)+\\ 
&&G^{0, 3}_{3, 1}\bigg(\begin{matrix}
1, z,  z \\
0, -
\end{matrix}\ \bigg|\ \frac{D_K}{4\pi ^2m}\bigg)\bigg\},
\\
H(z):= &&(2\pi)^{2z-1}|D_K |^{\frac{1-2z}{2}}
\sum_{m=1}^\infty\sigma'_{1-u-v, v_K}(m)
m^{z-1}\int_0^{\frac{4\pi^2m}{|D_K|}} t^{-z}J_0(2t^{1/2})\rm dt,\\
T(u, v):=&&\dfrac{h_K}{w_K |D_K|^{\frac{1}{2}}}\bigg(\dfrac{1}{u-1}+\dfrac{1}{v-1}\bigg)\z_K(u+v-1),
\end{eqnarray*}
where $J_0$ denotes the Bessel function of first kind.
The main result on the product of two Dedekind zeta functions is the following: 
\begin{theorem}\label{mainthm}
Let $\R (u)>-1, ~ \R (v)>-1,~ \R (u+v)>0,~ u\ne 1, ~ v\ne 1$ and $u+v\ne 2$. Then for any quadratic number field $K$,
\begin{numcases}{
\z_K(u)\z_K(v)=}4R_KT(u,v)+Z(u)+Z(v), \text{ if } K \text{ is real}, \label{realquadratic} \\
2\pi T(u, v)-H(u)-H(v), \text{ if } K \text{ is imaginary}. \label{imaginaryquadratic}
\end{numcases}
\end{theorem}
\noindent
\textbf{Remarks.} Here, we make a few remarks before proceeding further.
\begin{itemize}
\item[(I)] For $K=\mathbb{Q}$, Theorem \ref{mainthm} gives a formula for the product of two Riemann zeta functions, which is slightly different from Theorem \ref{Wilton}. 
\item[(II)] A similar line of argument as adopted in the proof should entail a version of Theorem \ref{mainthm} for higher degree number fields using other Meijer G-functions as well as  
for the product of finite number ($>2$) of Dedekind zeta functions. Obviously, the complications will increase manifold.
\item[(III)] One can extend the range of Theorem \ref{mainthm} using analytic continuation on both sides. 
\end{itemize}

\subsection{Special values of Dedekind zeta functions: Earlier works}\hspace{6cm}$   $
The values of Dedekind zeta functions at rational integers are closely related to the algebraic characters of the attached number fields. These values are linked in mysterious ways to the underlying geometry and often seem to dictate the most important properties of the objects to which the Dedekind zeta functions are connected. To represent these values as clearly as possible is one of the most important tasks in the present years. Over the decades, many mathematicians have been working on this project. 

In 1734, L. Euler established  a celebrated result, which states that
$\zeta(2m)$
is a rational multiple of $\pi^{2m}$ for all natural numbers $m$. As a generalization of this result for any totally real number field $K$ with discriminant $D_K$, C. L. Siegel established the following in \cite{SI1937}: 
$$\z_K(2m)=\text{rational}\times \frac{\pi^{2dm}}{\sqrt{D_K}},$$
where $d$ is the degree of $K$ over $\mathbb{Q}$. 
C. L. Siegel developed (using finite dimensionality of elliptic modular forms) an innovative method in \cite{SI1969} for computing the special values of Dedekind zeta function attached to a totally real number field at negative odd integers, and a particular interesting case is at $-1$. However, evaluation of special values of a Dedekind zeta function by means of this method requires complicated computations.
The residue of a Dedekind zeta function attached to an arbitrary number field $K$ at $1$ gives the analytic class number formula that involves $h_K$, $R_K$, $w_K$ and the absolute value of $D_K$. The problem of expressing the special values of Dedekind zeta functions in terms of elementary functions was first studied by E. Hecke in \cite{HE1921}. He obtained an elementary expression involving Dedekind sum for the relative class number of totally complex quadratic extensions of a real quadratic number field. C. Meyer \cite{ME1957, ME1966} and his student H. Lange \cite{LA1968} evaluated $L$-functions for real quadratic number fields in terms of Dedekind sums. The precise history of the special values of zeta and $L$-functions attached to a number field is quite lengthy and partly complicated due to the great popularity this subject is enjoying. D. Zagier is one of them who have made most important contributions in this subject.
In \cite{ZA1975}, he obtained an expression for the values of Dedekind zeta functions attached to real quadratic number fields at negative integers using Kronecker limit formula. Siegel's method in \cite{SI1969} has been exploited by  D. Zagier to give an elementary expression for $\z_K(1-2s)$ in \cite{ZA1976}, where $K$ is a real quadratic number field and $s$ is a positive integer, which involves only rational integers, but neither algebraic numbers nor norms of ideals. In 1976, T. Shintani \cite{SH1976}, using an impressive linear programming method, expressed Dedekind zeta function attached to a totally real number field into a finite sum of functions, each given by a Dirichlet series whose meromorphic continuation assumes rational values at negative integers. He also obtained a formula to express $\z_K(-n)$ for $ n = 0, 1, 2, \cdots$. However, little is known about the number $\z_K(2m)$ for the number field $K$ which is not totally real. In 1986, D. Zagier \cite{ZA1986} obtained a formula for $\z_K(2)$ for an arbitrary number field $K$. His method is completely geometric, where it involves the interpretation of $\z_K(2)$ as the volume of a hyperbolic manifold. Since it is the only $\z_K(2)$ which can be interpreted geometrically, this method fails to obtain the formula for $\z_K(2m)$ as well as for $\z_K(2m+1)$ for all natural numbers $m$. 
In the same paper, D. Zagier found a nice expression for $\z_K(2)$ attached to an imaginary quadratic number field $K$. 

In this article, we find some expressions for the values of Dedekind zeta functions attached to arbitrary real (resp. imaginary) quadratic number fields at any positive integers. Zagier's result on $\zeta_K(2)$ plays a crucial role in finding these expressions.

\subsection{Structure of the paper} \hspace{9cm}$   $
This paper comprises of $5$ sections. In \S 2, we discuss the main ingredients that are needed to prove our results. In \S 3, we provide the proof of Theorem \ref{mainthm}. In \S 4, we discuss a work of D. Zagier on the values of Dedekind zeta functions. We also discuss the expressions for the values of Dedekind zeta function for both real and imaginary quadratic number fields at any positive integers. In \S 5, we discuss some important remarks.  

\section{Background setup}
Riesz means were introduced by M. Riesz \cite{HR64} and have been studied in connection with summability of Fourier series and that of Dirichlet series (for details, we refer \cite{CS16, CM52}). For a given increasing sequence $\{\lambda_n\}$ of real numbers and a given sequence $\{\alpha_n\}$ of complex numbers, the Riesz sum of order $\ka$ is defined by
\begin{eqnarray*}
\mathcal{A}^\ka(x)=\mathcal{A}_\lambda^\ka (x)&=& \sideset{}{'}\sum_{\lambda_n\leq x}(x-\lambda_n)^\ka\alpha_n\\
&=&\ka\int_0^x(x-t)^{\ka-1}\mathcal{A}_\lambda (t){\rm d}t\\
&=&\int_0^x(x-t)^{\ka}{\rm d}\mathcal{A}_\lambda (t)
\end{eqnarray*}
with $\mathcal{A}_\lambda(x)=\mathcal{A}_\lambda^0(x)=\sideset{}{'}\sum\limits_{\lambda_n\leq x}\alpha_n$, where the prime on the summation sign means that when $\lambda_n=x$, the corresponding term is to be halved.
Sometimes normalized $\frac{1}{\Gamma(\ka+1)}\mathcal{A}^\ka(x)$ that appears in the Perron's formula 
$$
\frac{1}{\Gamma(\ka+1)}\sideset{}{'}\sum\limits_{\lambda_n\leq x}\alpha_n(x-\lambda_n)^\ka=\frac{1}{2\pi i}\int_{(C)}\frac{\Gamma(w)\varphi(w)x^{\ka+w}}{\Gamma(w+\ka+1)}{\rm d}w
$$
is also called the Riesz sum of order $\ka$ with $\varphi(w)=\sum\limits^\infty_{n=1}\frac{\alpha_n}{\lambda_n^w}$. 

We consider the Dirichlet series 
$$
\varphi(s)=\sum\limits^\infty_{n=1}\frac{\alpha_n}{\lambda_n^s}~~ \text{and}~~ \Phi(s)=\sum\limits^\infty_{m=1}\frac{a_m}{\gamma_m^s},
$$
where $\{\lambda_n\}$ and $\{\gamma_m\}$ are increasing sequences of real numbers. Here $\alpha_n$ and $a_m$ are complex numbers, and we assume that these series are continued to meromorphic functions over the whole complex plane, and that they satisfy suitable growth conditions. 

Further we consider the following integral,
\begin{eqnarray*}
\F_{(C)}^\ka (\varphi(u), \Phi(v);x)&=&\frac{1}{2\pi i}\int_{(C)}\frac{\Gamma(w)}{\Gamma(w+\ka+1)}\varphi(u+w)\Phi(v-w)x^{w+\ka}{\rm d}w\\
&=& \frac{1}{2\pi i}\sum_{m=1}^\infty \frac{a_m}{\gamma_m^v}\int_{(C)}\frac{\Gamma(w)}{\Gamma(w+\ka+1)}\sum_{n=1}^\infty\frac{\alpha_n}{\lambda_n^{u+w}} \gamma_m^wx^{w+\ka}{\rm d}w\\
&=& \frac{1}{2\pi i}\sum_{m=1}^\infty \frac{a_m}{\gamma_m^{v+\ka}}\int_{(C)}\frac{\Gamma(w)}{\Gamma(w+\ka+1)}\sum_{n=1}^\infty\frac{\alpha_n}{\lambda_n^{u+w}} (\gamma_m x)^{w+\ka}{\rm d}w\\
&=& \frac{1}{\Gamma(\ka+1)}\sum_{m=1}^\infty \frac{a_m}{\gamma_m^{v+\ka}}\sideset{}{'}\sum_{\lambda_n\leq \gamma_mx}\alpha_n\lambda_n^{-u}(\gamma_mx-\lambda_n)^\ka.
\end{eqnarray*}
Similarly we have,
\begin{eqnarray*}
\F_{(C)}^\ka (\Phi(v), \varphi(u);x)&=&\frac{1}{\Gamma(\ka+1)}\sum_{n=1}^\infty \frac{\alpha_n}{\lambda_n^{u+\ka}}\sideset{}{'}\sum_{\gamma_m\leq \lambda_nx}a_m\gamma_m^{-v}(\lambda_nx-\gamma_m)^\ka.
\end{eqnarray*}
Let $\varphi(u)=\z_K(u)$ and $\Phi(v)=\z_K(v)$ with order $\ka=1$. Then by using the above expressions we obtain,
\begin{eqnarray*}
& &\F_{(C)}^1 (\z_K(u), \z_K(v);x)+\F_{(C)}^1 (\z_K(v), \z_K(u);x)\\
&=&\sum_{m=1}^\infty \frac{v_K(m)}{m^{v+1}}\sideset{}{'}\sum_{n\leq mx}\frac{v_K(n)}{n^u}(mx-n)+\sum_{n=1}^\infty \frac{v_K(n)}{n^{u+1}}\sideset{}{'}\sum_{m\leq nx}\frac{v_K(m)}{m^v}(nx-m).
\end{eqnarray*}
We now differentiate with respect to $x$ and then substitute $x=1$ to get 
\begin{eqnarray}\label{split1}
& &(\F_{(C)}^1 (\z_K(u), \z_K(v);1))'+(\F_{(C)}^1 (\z_K(v), \z_K(u);1))'\nonumber\\
&=&\sum_{m=1}^\infty \frac{v_K(m)}{m^{v}}\sideset{}{'}\sum_{n\leq m}\frac{v_K(n)}{n^u}+\sum_{n=1}^\infty \frac{v_K(n)}{n^{u}}\sideset{}{'}\sum_{m\leq n}\frac{v_K(m)}{m^v}.\nonumber
\end{eqnarray}
M. Nakajima \cite{NA03} introduced a dissection (termed as `Nakajima dissection') involving the splitting of a double sum 
$$
\sum\limits_{m, n=1}^\infty=\sum_{m=1}^\infty \sideset{}{'}\sum_{n\leq m}+\sum_{n=1}^\infty\sideset{}{'}\sum_{m\leq n}.
$$
This dissection technique immediately gives,
\begin{eqnarray}\label{split2}
&&(\F_{(C)}^1 (\varphi(u), \Phi(v);1))'+(\F_{(C)}^1 (\Phi(v), \varphi(u);1))'\nonumber\\
&=&\sum_{m=1}^\infty \frac{v_K(m)}{m^{v}}\sum_{n=1}^\infty\frac{v_K(n)}{n^u}\nonumber\\
&=&\z_K(u)\z_K(v).
\end{eqnarray}
\section{Proof of the Theorem 1.1}
The Riesz sum set up, Perron's formula and Nakajima dissection are the main ingredients in this proof.
Assume that 
\begin{equation}\label{eq2.1}
 \F^\ka_{ (-B)}(\z_{K}(u), \z_{K}(v); x)=\frac{1}{2\pi i}\int_{(-B)} \frac{\Gamma (w)}{\Gamma (w+\ka+1)}\z_K(u+w)\z_K(v-w)x^{w+\ka}{\rm d}w
\end{equation}
with $\max \{ 0, \R(u)-1\} < B <\frac{3}{2}+\ka$ and $0<\R(v)-1<B<\frac{3}{2}+\ka$.
We assume additionally that $B < \R(u) - 1/2$ and $B < \R(v) - 1/2$. 
The functional equation satisfied by $\z_{K}(u+w)$ is 
\begin{eqnarray}\label{eq2.2}
\z_{K}(u+w)&=&2^{-r_2(1-2u-2w)}\pi^{-\frac{d}{2}(1-2u-2w)}|D_K|^{\frac{1-2u-2w}{2}}\frac{\Gamma(\frac{1-u-w}{2})^{r_1}}{\Gamma(\frac{u+w}{2})^{r_1}} \nonumber \\
&\times& \frac{\Gamma(1-u-w)^{r_2}}{\Gamma(u+w)^{r_2}} \z_K(1-u-w).
\end{eqnarray}
Now we have,
\begin{eqnarray}\label{eq2.3}
\frac{\Gamma(\frac{1-u-w}{2})^{r_1}}{\Gamma(\frac{u+w}{2})^{r_1}}&=&\frac{\big[ \Gamma(\frac{1-u-w}{2})\Gamma(1-\frac{u+w}{2})\big]^{r_1}}{\big[ \Gamma(\frac{u+w}{2})\Gamma(1-\frac{u+w}{2})\big]^{r_1}}\nonumber\\
&=&\frac{ \big[ \Gamma(\frac{1-u-w}{2})\Gamma(\frac{1}{2}+\frac{1-u-w}{2})\big]^{r_1}}{\bigg(\frac{\pi}{\sin{\frac{\pi}{2}(u+w)}}\bigg)^{r_1}}\nonumber \\
&=& \bigg( \frac{2^{u+w}}{\sqrt{\pi}}\bigg)^{r_1}\bigg(\sin\frac{\pi}{2}(u+w)\bigg)^{r_1}\Gamma(1-u-w)^{r_1}.
\end{eqnarray}
Similarly,
\begin{eqnarray}\label{eq2.4}
\frac{\Gamma(1-u-w)^{r_2}}{\Gamma(u+w)^{r_2}}&=& \frac{\Gamma(1-u-w)^{2r_2}}{\{\Gamma(u+w)\Gamma(1-u-w)\}^{r_2}}\nonumber\\
&=&\bigg(\frac{\sin \pi (u+w)}{\pi}\bigg)^{r_2}\Gamma(1-u-w)^{2r_2}.
\end{eqnarray}
Using (\ref{eq2.3}) and (\ref{eq2.4}) in (\ref{eq2.2}), we obtain 
\begin{eqnarray}
\z_K(u+w)&=& 2^{d(u+w)}\pi^{-d(1-u-w)}|D_K|^{\frac{1-2u-2w}{2}}\{\sin\frac{\pi}{2}(u+w)\}^{r_1+r_2}\nonumber \\
&\times& \{\cos\frac{\pi}{2}(u+w)\}^{r_2}\Gamma(1-u-w)^d \z_K(1-u-w). \nonumber
\end{eqnarray}
Thus (\ref{eq2.1}) becomes
\begin{eqnarray}\label{eq2.5}
\F^\ka_{ (-B)}(\z_{K}(u), \z_{K}(v); x)&=&\frac{1}{2\pi i}\int_{(-B)}S_\ka(w)f(w)2^{d(u+w)}\pi^{-d(1-u-w)}\nonumber \\
& \times& |D_K|^{\frac{1-2u-2w}{2}}\{\sin\frac{\pi}{2}(u+w)\}^{r_1+r_2}\{\cos\frac{\pi}{2}(u+w)\}^{r_2}\nonumber \\
&\times & \Gamma(1-u-w)^dx^{w+\ka}{\rm d}w.
\end{eqnarray}
Here $S_\ka(w)=\frac{\Gamma (w)}{\Gamma (w+\ka+1)}$ and $f(w)=\z_K(1-u-w)\z_K(v-w)$.

We now consider $K$ to be a quadratic number field. Then (\ref{eq2.5}) becomes
\begin{eqnarray}
\F^\ka_{ (-B)}(\z_{K}(u), \z_{K}(v); x)&=&\frac{1}{2\pi i}\int_{(-B)}S_\ka(w)f(w)2^{2(u+w)}\pi^{-2(1-u-w)}\nonumber \\
&\times & |D_K|^{\frac{1-2u-2w}{2}}\{\sin\frac{\pi}{2}(u+w)\}^{r_1+r_2}\{\cos\frac{\pi}{2}(u+w)\}^{r_2}\nonumber \\
&\times & \Gamma(1-u-w)^2x^{w+\ka}{\rm d}w.\nonumber
\end{eqnarray}
Employing change of variables we have
\begin{eqnarray}\label{eq2.6}
\F^\ka_{ (-B)}(\z_{K}(u), \z_{K}(v); x)&=&\frac{2}{\pi i}\int_{(B)}S_\ka(-z)f(-z)(2\pi)^{-2(1-u+z)}\nonumber \\
& \times & |D_K|^{\frac{1-2u+2z}{2}}\{\sin\frac{\pi}{2}(u-z)\}^{r_1+r_2}\{\cos\frac{\pi}{2}(u-z)\}^{r_2}\nonumber \\
&\times & \Gamma(1-u+z)^2x^{-z+\ka}{\rm d}z
\end{eqnarray}
where 
\begin{eqnarray}
f(-z)&=&\z_K(1-u+z)\z_K(v+z)\nonumber\\
&=& \sum_{m=1}^\infty \frac{v_K (m)}{m^{1-u+z}}\sum_{m=1}^\infty \frac{v_K (m)}{m^{v+z}}\nonumber\\
&=& \sum_{m=1}^\infty \frac{\sum_{d\mid m}d^{-v}\frac{m^{u-1}}{d^{u-1}}v_K(d)v_K(\frac{m}{d})}{m^z}\nonumber\\
&=& \sum_{m=1}^\infty \frac{\sum_{d\mid m}d^{1-u-v}v_K(d)v_K(\frac{m}{d})}{m^{1-u+z}}\nonumber\\
&=&\sum_{m=1}^\infty \frac{\sigma'_{1-u-v, v_K}(m)}{m^{1-u+z}}.\nonumber
\end{eqnarray}
Therefore (\ref{eq2.6}) implies that
\begin{eqnarray}
\F^\ka_{ (-B)}(\z_{K}(u), \z_{K}(v); x)&=&\frac{2}{\pi i} (2\pi)^{2(u-1)}|D_K |^{\frac{1-2u}{2}}\sum_{m=1}^\infty \sigma'_{1-u-v, v_K}(m)m^{u-1}\nonumber \\
& \times& \int_{(B)}S_\ka(-z)\bigg(\frac{|D_K|}{4\pi^2 m}\bigg)^{z}\{\sin\frac{\pi}{2}(u-z)\}^{r_1+r_2}
\nonumber \\
&\times & 
\{\cos\frac{\pi}{2}(u-z)\}^{r_2} \Gamma(1-u+z)^2x^{-z+\ka}{\rm d}z.\nonumber
\end{eqnarray}
Putting $\ka=1$, we have 
\begin{eqnarray}\label{eq2.7}
\F^1_{ (-B)}(\z_{K}(u), \z_{K}(v); x)&=&\frac{2}{\pi i} (2\pi)^{2(u-1)}|D_K |^{\frac{1-2u}{2}}\sum_{m=1}^\infty \sigma'_{1-u-v, v_K}(m)m^{u-1}\nonumber \\
&\times & \int_{(B)}\bigg(\frac{ |D_K |}{4\pi^2 m}\bigg)^{z}\{\sin\frac{\pi}{2}(u-z)\}^{r_1+r_2}
\nonumber \\
& \times& 
\{\cos\frac{\pi}{2}(u-z)\}^{r_2} \Gamma(1-u+z)^2\frac{x^{-z+1}}{z(z-1)}{\rm d}z.\nonumber\\
\end{eqnarray}
We treat real and imaginary cases separately.

{\bf Case-I}: {\underline{Real quadratic number field:}}

In this case $r_1=2$ and $r_2=0$. Thus (\ref{eq2.7}) gives
\begin{eqnarray}\label{eq2.8}
\F^1_{ (-B)}(\z_{K}(u), \z_{K}(v); x)&=&-\frac{1}{2\pi i} (2\pi)^{2(u-1)}D_K^{\frac{1-2u}{2}}\sum_{m=1}^\infty \sigma'_{1-u-v, v_K}(m)m^{u-1}\nonumber \\
& \times& \int_{(B)}\bigg(\frac{ D_K }{4\pi^2m}\bigg)^{z}\bigg(e^{i\pi(u-z)}+e^{-i\pi(u-z)}-2\bigg)
\nonumber \\
& \times& 
\Gamma(1-u+z)^2\frac{x^{-z+1}}{z(z-1)}{\rm d}z.
\end{eqnarray}

Let 
$$
H_B(u;x)=\int_{(B)}\bigg(\frac{ D_K}{4\pi^2m}\bigg)^{z}e^{i\pi(u-z)} \Gamma(1-u+z)^2\frac{x^{-z+1}}{z(z-1)}{\rm d}z.
$$
We differentiate $H_B(u;x)$ with respect to $x$ and get,
\begin{eqnarray}
H'_B(u;x)&=& -e^{i\pi u}\int_{(B)}\bigg(\frac{ D_K}{4\pi^2mx}e^{-i\pi}\bigg)^{z} \frac{\Gamma(1-u+z)^2}{z}{\rm d}z\nonumber\\
&=& -e^{i\pi u}\int_{(B)}\bigg(\frac{ D_K}{4\pi^2mx}e^{-i\pi}\bigg)^{z} \frac{\Gamma(z)\Gamma(1-u+z)\Gamma(1-u+z)}{\Gamma(z+1)}{\rm d}z.\nonumber\\
\end{eqnarray}
Also let,
$$
\mathcal{H}_B(u;x)=\int_{(B)}\bigg(\frac{ D_K }{4\pi^2m}\bigg)^{z}e^{-i\pi(u-z)} \Gamma(1-u+z)^2\frac{x^{-z+1}}{z(z-1)}{\rm d}z.
$$
As before differentiating with respect to $x$ gives
\begin{eqnarray}
\mathcal{H}'_B(u;x)&=& -e^{-i\pi u}\int_{(B)}\bigg(\frac{ D_K}{4\pi^2mx}e^{i\pi}\bigg)^{z} \frac{\Gamma(1-u+z)^2}{z}{\rm d}z\nonumber\\
&=& -e^{-i\pi u}\int_{(B)}\bigg(\frac{ D_K}{4\pi^2mx}e^{i\pi}\bigg)^{z} \frac{\Gamma(z)\Gamma(1-u+z)\Gamma(1-u+z)}{\Gamma(z+1)}{\rm d}z.\nonumber\\
\end{eqnarray}
Further, let
$$
\mathfrak{H}_B(u;x)=2\int_{(B)}\bigg(\frac{ D_K}{4\pi^2m}\bigg)^{z} \Gamma(1-u+z)^2\frac{x^{-z+1}}{z(z-1)}{\rm d}z.
$$
As before,
\begin{eqnarray}
\mathfrak{H}'_B(u;x)&=& -2\int_{(B)}\bigg(\frac{D_K}{4\pi^2mx}\bigg)^{z} \frac{\Gamma(1-u+z)^2}{z}{\rm d}z\nonumber\\
&=& -2\int_{(B)}\bigg(\frac{ D_K }{4\pi^2mx}\bigg)^{z} \frac{\Gamma(z)\Gamma(1-u+z)\Gamma(1-u+z)}{\Gamma(z+1)}{\rm d}z.\nonumber\\
\end{eqnarray}
The convergence of the integrals of $H'_B(u;x)$, $\mathcal{H}'_B(u;x)$ and $\mathfrak{H}'_B(u;x)$  follows for $$ \max \{0, \R(u) - 1\} < B < \R(u) - 1/2$$ which we have assumed in the beginning of the proof of Theorem $1.1$. Now applying the definition of Meijer G-function and summing the integrals we have

\begin{eqnarray}
& &H'_B(u;x)+\mathcal{H}'_B(u;x)+\mathfrak{H}'_B(u;x)=-2\pi i e^{i\pi u}G^{0, 3}_{3, 1}\bigg(\begin{matrix}
1, u,  u \\
0, -
\end{matrix}\ \bigg|\ \frac{D_K}{4\pi ^2mx}e^{-i\pi}\bigg)\nonumber \\ 
&-&2\pi i e^{-i\pi u}G^{0, 3}_{3, 1}\bigg(\begin{matrix}
1, u,  u \\
0, -
\end{matrix}\ \bigg|\ \frac{D_K}{4\pi ^2mx}e^{i\pi}\bigg)-4\pi i G^{0, 3}_{3, 1}\bigg(\begin{matrix}
1, u,  u \\
0, -
\end{matrix}\ \bigg|\ \frac{D_K}{4\pi ^2mx}\bigg).
\end{eqnarray}
Therefore (\ref{eq2.8}) becomes 
\begin{eqnarray}\label{eq2.9}
& &(\F^1_{ (-B)}(\z_{K}(u), \z_{K}(v); x))'=2(2\pi)^{2(u-1)}D_K ^{\frac{1-2u}{2}}\sum_{m=1}^\infty \sigma'_{1-u-v, v_K}(m)m^{u-1}\nonumber\\ 
& &\times
\bigg\{- e^{i\pi u}G^{0, 3}_{3, 1}\bigg(\begin{matrix}
1, u,  u \\
0, -
\end{matrix}\ \bigg|\ \frac{D_K}{4\pi ^2mx}e^{-i\pi}\bigg)- e^{-i\pi u}G^{0, 3}_{3, 1}\bigg(\begin{matrix}
1, u,  u \\
0, -
\end{matrix}\ \bigg|\ \frac{D_K}{4\pi ^2mx}e^{i\pi}\bigg)\nonumber\\ 
& &- 2 G^{0, 3}_{3, 1}\bigg(\begin{matrix}
1, u,  u \\
0, -
\end{matrix}\ \bigg|\ \frac{D_K}{4\pi ^2mx}\bigg)\bigg\}.
\end{eqnarray}
Now using Cauchy's residue theorem, we have,
\begin{eqnarray}
\F^1_{ (C)}(\z_{K}(u), \z_{K}(v); x)&=&\F^1_{ (-B)}(\z_{K}(u), \z_{K}(v); x) +x\z_K(u)\z_K(v)- \nonumber \\ 
& & \z_K(u-1) \z_K(v+1)+\frac{\z_K(u+v-1)x^{2-u}}{(u-2)(u-1)}\frac{4h_K R_K}{w_K D_K^{\frac{1}{2}}}.\nonumber 
\end{eqnarray}
Further differentiating with respect to $x$ gives,
\begin{eqnarray}\label{eq2.10}
(\F^1_{ (C)}(\z_{K}(u), \z_{K}(v); x))'&=& (\F^1_{ (-B)}(\z_{K}(u), \z_{K}(v); x))' +\z_K(u)\z_K(v)-\nonumber \\
& & \frac{\z_K(u+v-1)x^{1-u}}{u-1}\frac{4h_K R_K}{w_K D_K^{\frac{1}{2}}}.
\end{eqnarray}
Now (\ref{eq2.9}) and (\ref{eq2.10}) together imply,  
\begin{eqnarray}
& &(\F^1_{ (C)}(\z_{K}(u), \z_{K}(v); x))'=2(2\pi)^{2(u-1)}D_K^{\frac{1-2u}{2}}\sum_{m=1}^\infty \sigma'_{1-u-v, v_K}(m)m^{u-1}\nonumber \\
& & \times \bigg\{- e^{i\pi u}G^{0, 3}_{3, 1}\bigg(\begin{matrix}
1, u,  u \\
0, -
\end{matrix}\ \bigg|\ \frac{D_K}{4\pi ^2mx}e^{-i\pi}\bigg)- e^{-i\pi u}G^{0, 3}_{3, 1}\bigg(\begin{matrix}
1, u,  u \\
0, -
\end{matrix}\ \bigg|\ \frac{D_K}{4\pi ^2mx}e^{i\pi}\bigg)\nonumber\\ 
& &- 2 G^{0, 3}_{3, 1}\bigg(\begin{matrix}
1, u,  u \\
0, -
\end{matrix}\ \bigg|\ \frac{D_K}{4\pi ^2mx}\bigg)\bigg\} +\z_K(u)\z_K(v)-
\frac{\z_K(u+v-1)x^{1-u}}{u-1}\nonumber \\
& &\times \frac{4h_K R_K}{w_K D_K^{\frac{1}{2}}}.\nonumber
\end{eqnarray}
We now substitute $x=1$, 
 \begin{eqnarray}\label{eq2.13}
& &(\F^1_{ (C)}(\z_{K}(u), \z_{K}(v); 1))'=2(2\pi)^{2(u-1)}D_K ^{\frac{1-2u}{2}}\sum_{m=1}^\infty \sigma'_{1-u-v, v_K}(m)m^{u-1}\nonumber \\
& & \times \bigg\{- e^{i\pi u}G^{0, 3}_{3, 1}\bigg(\begin{matrix}
1, u,  u \\
0, -
\end{matrix}\ \bigg|\ \frac{D_K}{4\pi ^2m}e^{-i\pi}\bigg)- e^{-i\pi u}G^{0, 3}_{3, 1}\bigg(\begin{matrix}
1, u,  u \\
0, -
\end{matrix}\ \bigg|\ \frac{D_K}{4\pi ^2m}e^{i\pi}\bigg)\nonumber\\ 
& &- 2 G^{0, 3}_{3, 1}\bigg(\begin{matrix}
1, u,  u \\
0, -
\end{matrix}\ \bigg|\ \frac{D_K}{4\pi ^2m}\bigg)\bigg\} +\z_K(u)\z_K(v)-
\frac{\z_K(u+v-1)}{u-1}\nonumber \\
& &\times \frac{4h_K R_K}{w_K D_K^{\frac{1}{2}}}.
\end{eqnarray}
Similarly we have,
\begin{eqnarray}\label{eq2.14}
&&(\F^1_{ (C)}(\z_{K}(v), \z_{K}(u); 1))'= 2(2\pi)^{2(v-1)}D_K^{\frac{1-2v}{2}}\sum_{m=1}^\infty \sigma'_{1-u-v, v_K}(m)m^{v-1}\nonumber \\
& & \times \bigg\{- e^{i\pi v}G^{0, 3}_{3, 1}\bigg(\begin{matrix}
1, v,  v \\
0, -
\end{matrix}\ \bigg|\ \frac{D_K}{4\pi ^2m}e^{-i\pi}\bigg)- e^{-i\pi v}G^{0, 3}_{3, 1}\bigg(\begin{matrix}
1, v,  v \\
0, -
\end{matrix}\ \bigg|\ \frac{D_K}{4\pi ^2m}e^{i\pi}\bigg)\nonumber\\ 
& &- 2 G^{0, 3}_{3, 1}\bigg(\begin{matrix}
1, v,  v \\
0, -
\end{matrix}\ \bigg|\ \frac{D_K}{4\pi ^2m}\bigg)\bigg\} +\z_K(u)\z_K(v)-
\frac{\z_K(u+v-1)}{v-1}\nonumber \\
& &\times \frac{4h_K R_K}{w_K D_K^{\frac{1}{2}}}.
\end{eqnarray}
Now we add (\ref{eq2.13}) and (\ref{eq2.14}), and further utilize (\ref{split2}) to arrive at
\begin{eqnarray}
&&\z_K(u)\z_K(v)= 2\z_K(u)\z_K(v)-\frac{4h_K R_K}{w_K D_K^{\frac{1}{2}}}\bigg(\frac{1}{u-1}+\frac{1}{v-1}\bigg)\z_K(u+v-1)\nonumber\\
& & + 2(2\pi)^{2(u-1)}D_K ^{\frac{1-2u}{2}}\sum_{m=1}^\infty \sigma'_{1-u-v, v_K}(m)m^{u-1}\bigg\{- e^{i\pi u}G^{0, 3}_{3, 1}\bigg(\begin{matrix}
1, u,  u \\
0, -
\end{matrix}\ \bigg|\ \frac{D_K}{4\pi ^2m}e^{-i\pi}\bigg)\nonumber\\
&& - e^{-i\pi u}G^{0, 3}_{3, 1}\bigg(\begin{matrix}
1, u,  u \\
0, -
\end{matrix}\ \bigg|\ \frac{D_K}{4\pi ^2m}e^{i\pi}\bigg)- 2 G^{0, 3}_{3, 1}\bigg(\begin{matrix}
1, u,  u \\
0, -
\end{matrix}\ \bigg|\ \frac{D_K}{4\pi ^2m}\bigg)\bigg\}+2(2\pi)^{2(v-1)}\nonumber\\
&&\times  D_K ^{\frac{1-2v}{2}}\sum_{m=1}^\infty \sigma'_{1-u-v, v_K}(m)m^{v-1}\bigg\{- e^{i\pi v}G^{0, 3}_{3, 1}\bigg(\begin{matrix}
1, v,  v \\
0, -
\end{matrix}\ \bigg|\ \frac{D_K}{4\pi ^2m}e^{-i\pi}\bigg)\nonumber\\
&&- e^{-i\pi v}G^{0, 3}_{3, 1}\bigg(\begin{matrix}
1, v,  v \\
0, -
\end{matrix}\ \bigg|\ \frac{D_K}{4\pi ^2m}e^{i\pi}\bigg)- 2 G^{0, 3}_{3, 1}\bigg(\begin{matrix}
1, v,  v \\
0, -
\end{matrix}\ \bigg|\ \frac{D_K}{4\pi ^2m}\bigg)\bigg\}.\nonumber
\end{eqnarray}
Which implies that,
\begin{eqnarray}
\z_K(u)\z_K(v)&=&\frac{4h_K R_K}{w_K D_K^{\frac{1}{2}}}\bigg(\frac{1}{u-1}+\frac{1}{v-1}\bigg)\z_K(u+v-1) +2(2\pi)^{2(u-1)}\nonumber\\
&\times& D_K^{\frac{1-2u}{2}}\sum_{m=1}^\infty \sigma'_{1-u-v, v_K}(m)m^{u-1}\bigg\{ e^{i\pi u}G^{0, 3}_{3, 1}\bigg(\begin{matrix}
1, u,  u \\
0, -
\end{matrix}\ \bigg|\ \frac{D_K}{4\pi ^2m}e^{-i\pi}\bigg)\nonumber\\
&+& e^{-i\pi u}G^{0, 3}_{3, 1}\bigg(\begin{matrix}
1, u,  u \\
0, -
\end{matrix}\ \bigg|\ \frac{D_K}{4\pi ^2m}e^{i\pi}\bigg)+ 2 G^{0, 3}_{3, 1}\bigg(\begin{matrix}
1, u,  u \\
0, -
\end{matrix}\ \bigg|\ \frac{D_K}{4\pi ^2m}\bigg)\bigg\}+2(2\pi)^{2(v-1)}\nonumber\\
&\times&  D_K ^{\frac{1-2v}{2}}\sum_{m=1}^\infty \sigma'_{1-u-v, v_K}(m)m^{v-1}\bigg\{ e^{i\pi v}G^{0, 3}_{3, 1}\bigg(\begin{matrix}
1, v,  v \\
0, -
\end{matrix}\ \bigg|\ \frac{D_K}{4\pi ^2m}e^{-i\pi}\bigg)\nonumber\\
&+& e^{-i\pi v}G^{0, 3}_{3, 1}\bigg(\begin{matrix}
1, v,  v \\
0, -
\end{matrix}\ \bigg|\ \frac{D_K}{4\pi ^2m}e^{i\pi}\bigg)+2 G^{0, 3}_{3, 1}\bigg(\begin{matrix}
1, v,  v \\
0, -
\end{matrix}\ \bigg|\ \frac{D_K}{4\pi ^2m}\bigg)\bigg\}.\nonumber
\end{eqnarray}
This concludes the proof of this case.\vspace{3mm}\\
{\bf Case-II}: {\underline{Imaginary quadratic number field:}}\\
In this case $r_1=0$ and $r_2=1$. Thus (\ref{eq2.7}) implies that,
\begin{eqnarray}\label{eq2.15}
\F^1_{ (-B)}(\z_{K}(u), \z_{K}(v); x)&=&\frac{1}{\pi i} (2\pi)^{2(u-1)}|D_K |^{\frac{1-2u}{2}}\sum_{m=1}^\infty \sigma'_{1-u-v, v_K}(m)m^{u-1}\nonumber \\
&\times & \int_{(B)}\bigg(\frac{ |D_K |}{4\pi^2m}\bigg)^{z}\sin\pi(u-z)
\Gamma(1-u+z)^2 \frac{x^{-z+1}}{z(z-1)}{\rm d}z\nonumber\\
&=& -\frac{1}{2\pi} (2\pi)^{2(u-1)}|D_K |^{\frac{1-2u}{2}}\sum_{m=1}^\infty \sigma'_{1-u-v, v_K}(m)m^{u-1}\nonumber \\
&\times & (I_B(u;x)-\mathcal{I}_B(u;x)),
\end{eqnarray}
where,
$$
I_B(u;x)=\int_{(B)}\bigg(\frac{ |D_K |}{4\pi^2m}\bigg)^{z}e^{i\pi(u-z)} \Gamma(1-u+z)^2\frac{x^{-z+1}}{z(z-1)}{\rm d}z,
$$
and
$$
\mathcal{I}_B(u;x)=\int_{(B)}\bigg(\frac{ |D_K |}{4\pi^2m}\bigg)^{z}e^{-i\pi(u-z)} \Gamma(1-u+z)^2\frac{x^{-z+1}}{z(z-1)}{\rm d}z.
$$
Proceeding analogously we differentiate $I_B(u;x)$ with respect to $x$ and apply the definition of Meijer G-function to get
\begin{eqnarray}\label{eqa}
I'_B(u;x)&=& -e^{i\pi u}\int_{(B)}\bigg(\frac{ |D_K |}{4\pi^2mx}e^{-i\pi}\bigg)^{z} \frac{\Gamma(1-u+z)^2}{z}{\rm d}z\nonumber\\
&=& -e^{i\pi u}\int_{(B)}\bigg(\frac{ |D_K |}{4\pi^2mx}e^{-i\pi}\bigg)^{z} \frac{\Gamma(z)\Gamma(1-u+z)\Gamma(1-u+z)}{\Gamma(z+1)}{\rm d}z\nonumber\\
&=& -2\pi i e^{i\pi u}G^{0, 3}_{3, 1}\bigg(\begin{matrix}
1, u,  u \\
0, -
\end{matrix}\ \bigg|\ \frac{|D_K|}{4\pi ^2mx}e^{-i\pi}\bigg).
\end{eqnarray}
Again differentiating $\mathcal{I}_B(u;x)$ with respect to $x$ and applying the definition of Meijer G-function  we get
\begin{eqnarray}\label{eqb}
\mathcal{I}'_B(u;x)&=& -e^{-i\pi u}\int_{(B)}\bigg(\frac{|D_K |}{4\pi^2mx}e^{i\pi}\bigg)^{z} \frac{\Gamma(1-u+z)^2}{z}{\rm d}z\nonumber\\
&=& -2\pi i e^{-i\pi u}G^{0, 3}_{3, 1}\bigg(\begin{matrix}
1, u,  u \\
0, -
\end{matrix}\ \bigg|\ \frac{|D_K|}{4\pi ^2mx}e^{i\pi}\bigg).
\end{eqnarray}
Here we can apply the definition of Meijer G-function since the integrals in $I'_B(u;x)$ and $\mathcal{I}'_B(u;x)$ are convergent for 
$$ \max \{0, \R(u) - 1\} < B < \R(u) - 1/2$$
which we have assumed in the beginning of the proof of Theorem $1.1$.  
We now recall the following relation (see, relation (14) \cite[pp.210]{Erdelyi}) which is needed to proceed further:
\begin{eqnarray}\label{relation14}
G^{m, \ n}_{p, \ q}\bigg(\begin{matrix}
a_1, \cdots, a_p \\
b_1, \cdots, b_q
\end{matrix} \ \bigg|\ z\bigg)&=& \frac{1}{2\pi i}\bigg\{ e^{\pi i a_{n+1}}G^{m, \ n+1}_{p, \ q}\bigg(\begin{matrix}
a_1, \cdots, a_p \\
b_1, \cdots, b_q
\end{matrix} \ \bigg|\ z e^{-\pi i}\bigg)\nonumber\\
&-& e^{-\pi i a_{n+1}}G^{m, \ n+1}_{p, \ q}\bigg(\begin{matrix}
a_1, \cdots, a_p \\
b_1, \cdots, b_q
\end{matrix} \ \bigg|\ z e^{\pi i}\bigg)\bigg\} 
\end{eqnarray}
whenever $n\leq p-1$.

Subtracting \eqref{eqb} from \eqref{eqa} and then applying \eqref{relation14} we obtain,
\begin{eqnarray}\label{eqi-i}
I'_B(u;x)-\mathcal{I}'_B(u;x)&=&-2\pi i\bigg\{ e^{i\pi u}G^{0, 3}_{3, 1}\bigg(\begin{matrix}
1, u,  u \\
0, -
\end{matrix}\ \bigg|\ \frac{|D_K|}{4\pi ^2mx}e^{-i\pi}\bigg)\nonumber\\
&-& e^{-i\pi u}G^{0, 3}_{3, 1}\bigg(\begin{matrix}
1, u,  u \\
0, -
\end{matrix}\ \bigg|\ \frac{|D_K|}{4\pi ^2mx}e^{i\pi}\bigg)\bigg\}\nonumber\\
&=&4\pi^2 G^{0, 2}_{3, 1}\bigg(\begin{matrix}
1, u,  u \\
0, -
\end{matrix}\ \bigg|\ \frac{|D_K|}{4\pi ^2mx}\bigg).
\end{eqnarray}
We further recall the following important relation (see, relation (9) \cite[pp.209]{Erdelyi}) to rewrite \eqref{eqi-i} in suitable form: 
\begin{eqnarray}\label{relation9}
G^{m, \ n}_{p, \ q}\bigg(\begin{matrix}
a_1, \cdots, a_p \\
b_1, \cdots, b_q
\end{matrix} \ \bigg|\ \frac{1}{z}\bigg)&=& G^{n, \ m}_{q, \ p}\bigg(\begin{matrix}
1-b_1, \cdots, 1-b_q \\
1-a_1, \cdots, 1-a_p
\end{matrix} \ \bigg|\ z\bigg).
\end{eqnarray}

Applying relation \eqref{relation9} in \eqref{eqi-i}, we obtain
  \begin{eqnarray}
I'_B(u;x)-\mathcal{I}'_B(u;x)&=&4\pi ^2 G^{2, 0}_{1, 3}\bigg(\begin{matrix}
1 \\
0, 1-u, 1-u
\end{matrix}\ \bigg|\ \frac{4\pi ^2mx}{|D_K|}\bigg).
\end{eqnarray}
Putting $X=\frac{4\pi^2mx}{|D_K|}$ in \begin{eqnarray*}
G^{2, 0}_{1, 3}\bigg(\begin{matrix}
1 \\
0, 1-u, 1-u
\end{matrix}\ \bigg|\ \frac{4\pi ^2mx}{|D_K|}\bigg),
\end{eqnarray*}
we have
\begin{eqnarray}\label{eqc}
G^{2, 0}_{1, 3}\bigg(\begin{matrix}
1 \\
0, 1-u, 1-u
\end{matrix}\ \bigg|\ X\bigg)&=&-\frac{1}{2\pi i} \int_{(B)}\frac{\Gamma(1-u-s)}{\Gamma(u+s)}\frac{X^s}{s}\rm ds\nonumber\\
&=&-\frac{1}{2\pi i}\int_0^X \int_{(B)}\frac{\Gamma(1-u-s)}{\Gamma(u+s)}t^{s-1}\rm ds~\rm dt\nonumber\\
&=&-\int_0^X G^{1, 0}_{0, 2}\bigg(\begin{matrix}
- \\
-u, -u
\end{matrix}\ \bigg|\ t\bigg)\rm dt\nonumber\\
&=&-\int_0^X t^{-u}J_0(2t^{1/2})\rm dt. 
\end{eqnarray}

We now differentiate \eqref{eq2.15} with respect to $x$ to get 
\begin{eqnarray}
(\F^1_{ (-B)}(\z_{K}(u), \z_{K}(v); x))'&=&(2\pi)^{2u-1}|D_K |^{\frac{1-2u}{2}} \sum_{m=1}^\infty\sigma'_{1-u-v, v_K}(m)\nonumber\\
&\times &
m^{u-1}\int_0^X t^{-u}J_0(2t^{1/2})\rm dt.\nonumber
\end{eqnarray}
Using Cauchy's residue theorem, we get 
\begin{eqnarray}
(\F^1_{ (C)}(\z_{K}(u), \z_{K}(v); x))'&=&(2\pi)^{2u-1}|D_K |^{\frac{1-2u}{2}} \sum_{m=1}^\infty\sigma'_{1-u-v, v_K}(m)\nonumber\\
&\times &
m^{u-1}\int_0^X t^{-u}J_0(2t^{1/2}){\rm dt}+ 
\z_{K} (u)  \z_{K} (v)\nonumber \\
&- & 
\frac{\z_K(u+v-1)x^{1-u}}{u-1}\frac{2\pi h_K}{w_K |D_K|^{\frac{1}{2}}}.\nonumber
\end{eqnarray}
Further substituting $x=1$,
\begin{eqnarray}\label{eq2.16}
(\F^1_{ (C)}(\z_{K}(u), \z_{K}(v); 1))'&=&(2\pi)^{2u-1}|D_K |^{\frac{1-2u}{2}} \sum_{m=1}^\infty\sigma'_{1-u-v, v_K}(m)\nonumber\\
&\times &
m^{u-1}\int_0^{\frac{4\pi^2m}{|D_K|}} t^{-u}J_0(2t^{1/2}){\rm dt}+\z_K(u)\z_K(v)\nonumber \\
&- & 
\frac{\z_K(u+v-1)}{u-1}\frac{2\pi h_K}{w_K |D_K|^{\frac{1}{2}}}.
\end{eqnarray}
Similarly we obtain,
\begin{eqnarray}\label{eq2.17}
(\F^1_{ (C)}(\z_{K}(v), \z_{K}(u); 1))'&=& (2\pi)^{2v-1}|D_K |^{\frac{1-2v}{2}} \sum_{m=1}^\infty\sigma'_{1-u-v, v_K}(m)\nonumber\\
&\times &
m^{v-1}\int_0^{\frac{4\pi^2m}{|D_K|}} t^{-v}J_0(2t^{1/2}){\rm dt}+\z_K(u)\z_K(v)\nonumber \\
&- & 
\frac{\z_K(u+v-1)}{v-1}\frac{2\pi h_K}{w_K |D_K|^{\frac{1}{2}}}.
\end{eqnarray}
Adding (\ref{eq2.16}) and (\ref{eq2.17}), and then applying (\ref{split2}) we obtain,
\begin{eqnarray}
\z_K(u)\z_K(v)&=& 2\z_K(u)\z_K(v)-\frac{2\pi h_K}{w_K |D_K|^{\frac{1}{2}}}\bigg(\frac{1}{u-1}+\frac{1}{v-1}\bigg)\z_K(u+v-1)\nonumber\\
&+ & (2\pi)^{2u-1}|D_K |^{\frac{1-2u}{2}} \sum_{m=1}^\infty\sigma'_{1-u-v, v_K}(m)
m^{u-1}\int_0^{\frac{4\pi^2m}{|D_K|}} t^{-u}J_0(2t^{1/2})\rm dt\nonumber \\ &+&(2\pi)^{2v-1}|D_K |^{\frac{1-2v}{2}} \sum_{m=1}^\infty\sigma'_{1-u-v, v_K}(m)
m^{v-1}\int_0^{\frac{4\pi^2m}{|D_K|}} t^{-v}J_0(2t^{1/2})\rm dt.\nonumber
\end{eqnarray}
Finally this gives the required result:
\begin{eqnarray}
\z_K(u)\z_K(v)&=&\frac{2\pi h_K}{w_K |D_K|^{\frac{1}{2}}}\bigg(\frac{1}{u-1}+\frac{1}{v-1}\bigg)\z_K(u+v-1) -(2\pi)^{2u-1}|D_K |^{\frac{1-2u}{2}} \nonumber\\
&\times& \sum_{m=1}^\infty\sigma'_{1-u-v, v_K}(m)
m^{u-1}\int_0^{\frac{4\pi^2m}{|D_K|}} t^{-u}J_0(2t^{1/2}){\rm dt}-(2\pi)^{2v-1}|D_K |^{\frac{1-2v}{2}} \nonumber\\
&\times&\sum_{m=1}^\infty\sigma'_{1-u-v, v_K}(m)
m^{v-1}\int_0^{\frac{4\pi^2m}{|D_K|}} t^{-v}J_0(2t^{1/2})\rm dt.\nonumber
\end{eqnarray}

\section{Values of Dedekind zeta functions}
The special values of Dedekind zeta functions attached to an arbitrary quadratic number field is not completely known. In this section, we will discuss about  some formulae for Dedekind zeta functions attached to both real and imaginary quadratic number fields at any positive integers. Zagier's result on $\zeta_K(2)$ is used to get these formulae. 
To state his result, we need to introduce  the following real valued function
$$
A(x)=\int_0^x \frac{1}{1+t^2}\log \frac{4}{1+t^2}dt.
$$
The precise result of D. Zagier \cite{ZA1986} for $\z_K(2)$ is the following:
\begin{thma}\label{Zagier1}
Let $K$ be an arbitrary number field with discriminant $D_K$, and let $r_1$ and $r_2$ be its numbers of real and complex places, respectively. Then 
\begin{equation*}
\z_K(2)=\frac{\pi^{2r_1+2r_2}}{\sqrt{|D_K|}}\times \underset{\text{finite}}{\sum_\nu c_\nu} A(x_{\nu,1})\cdots A(x_{\nu,r_2}). 
\end{equation*}
Here  $c_\nu$'s are rational numbers and $x_{\nu, j}$'s are real algebraic numbers. 
\end{thma} 

This formula does not work for finding $\z_K(2m)$ with $m\geq 2$. However, he conjectured that in this case too, an analogous result holds. We need a more general real-valued function to state the conjecture which appeared in \cite{ZA1986}. For each natural number $m$, let $A_m$ be the real-valued function defined by 
$$
A_m(x)=\frac{2^{2m-1}}{(2m-1)!}\int_0^\infty\frac{t^{2m-1}}{x\sinh^2t+x^{-1}\cosh^2t}dt.
$$
It is noted that the function $A_m(x)$ agrees with the function $A(x)$ when $m=1$.

\begin{conjecture}[Zagier]\label{conj1}
Let $K$ be an arbitrary number field with discriminant $D_K$, and let $r_1$ and $r_2$ be its numbers of real and complex places, respectively. Then $\z_K(2m)$ equals $\frac{\pi^{2m(r_1+r_2)}}{\sqrt{|D_K|}}$ times a rational linear combination of the products of $r_2$ values of $A_m(x)$ at algebraic arguments.
\end{conjecture}
In \cite{ZA1986}, he proved this conjecture for certain abelian number fields. More precisely, he proved the following:
\begin{thma}\label{Zagier3} 
If $K$ is an abelian number field, then Conjecture \ref{conj1} is true. In fact, the argument $x$ can be chosen of the form $x=\cot\frac{\pi t}{N}$, where $N$ is the smallest natural number such that $K\subset \mathbb{Q}(e^{2\pi i/N})$.
\end{thma}
For $m=1$, and for an imaginary quadratic number field $K$, Theorem \ref{Zagier3} reduces to the following simplest form:
\begin{equation}\label{zetac1}
\z_K(2)=\frac{\pi^2}{6\sqrt{|D_K|}}\times \underset{0<n<|D_K|}\sum\bigg(\frac{D_K}{n}\bigg) A\bigg(\cot \frac{\pi n}{|D_K|}\bigg),
\end{equation} 
where the arguments of $A(x)$ for $\gcd (n, D_K)=1$, are of degree $\phi(|D_K|)$ or $\phi(|D_K|/2)$ over $\mathbb{Q}$, where $\phi(.)$ is Euler phi function.
In particular, we have the following:
\begin{itemize}
\item[(i)] For $D_K=-7$,
\begin{eqnarray*}
\z_{\mathbb{Q}(\sqrt{-7})}(2)&=&\frac{\pi^2}{3\sqrt{7}}\bigg\{A\bigg(\cot\frac{\pi}{7}\bigg)+ A\bigg(\cot\frac{2\pi}{7}\bigg)+A\bigg(\cot\frac{4\pi}{7}\bigg)\bigg\}\\
&\cong &1.89484144897.
\end{eqnarray*}
\item[(ii)] For $D_K=-3$,
\begin{eqnarray*}
\z_{\mathbb{Q}(\sqrt{-3})}(2)&=&\frac{\pi^2}{3\sqrt{3}}A\big(\frac{1}{\sqrt{3}}\big)\\
&\cong & 1.28519145388.
\end{eqnarray*}
\end{itemize}
In fact,  D. Zagier expressed $\z_K(2)$ for an imaginary quadratic number field $K$ in a much simpler form in the same paper. More precisely, he proved:
\begin{thma}\label{Zagier4}
Let $K=\mathbb{Q}(\sqrt{D_K})$ be an imaginary quadratic number field of discriminant $D_K$. Then 
\begin{equation}\label{zetac2}
\z_K(2)=\frac{\pi^2}{6|D_K|\sqrt{|D_K|}}\times \underset{\text{finite}}{\sum_\nu \eta_\nu} A\bigg(\frac{\lambda_\nu}{\sqrt{|D_K|}}\bigg), 
\end{equation}
for some finite collection of numbers $\eta_\nu\in \mathbb{Z}$ and $\lambda_\nu \in \mathbb{Q}$ with $\underset{\nu}{\prod}\big(\lambda_\nu+\sqrt{|D_K|}\big)^{\eta_\nu}\in \mathbb{Q}$.
\end{thma}

\subsection{Real quadratic number field case}

This section deals with an expression for the values of Dedekind zeta functions attached to real quadratic number fields. We note that an expression for the values of Dedekind zeta functions at even positive integers can be obtained from two prominent results. One of them was given by Euler which entails the expression of zeta function at even positive integers as
\begin{equation}\label{eqn6.3}
\zeta(2n) = \frac{(-1)^{n+1}B_{2n}(2\pi)^{2n}}{2 \, (2n)!},
\end{equation}  
where $B_{2n}$ is the $2n$-th Bernoulli number.

Let $\chi$ be any Dirichlet character with modulus $f$ such that $\chi(-1)=1$. Then we can derive the following relation from a result of H. -W. Leopoldt (See, relation (12) in \cite{Leopoldt}):

\begin{equation}\label{eqn6.4}
L(2n, \chi) = (-1)^{n+1}\frac{\tau(\chi)}{2} \left(\frac{2\pi}{f}\right)^{2n} \frac{B_{2n ,\overline{\chi}}}{(2n)!},
\end{equation}
where $\tau(\chi)$ is the Gauss sum associated to $\chi$ and $B_{2n, \chi}$ is the $2n$-th generalized Bernoulli number associated to $\chi$.

Let $K = \mathbb{Q}(\sqrt{D})$ be a quadratic number field with discriminant
$D_K$. The Kronecker symbol $(D_K /.)$ is a real primitive Dirichlet character with modulus $|D_k|$. We denote this Dirichlet character by $\chi_{D_K}$ and it's associated Dirichlet $L$-series would be
\begin{equation}
L(s,\chi_{D_K}) =\sum_{m=1}^\infty(D_K /m) m^{-s}.\nonumber
\end{equation}
This is convergent in the half-plane $\sigma>0$.
We have 
\begin{equation}\label{eqn6.6}
\zeta_K(s) = \zeta(s) L(s, \chi_{D_K}).
\end{equation}
We now apply \eqref{eqn6.3} and \eqref{eqn6.4} in \eqref{eqn6.6} to obtain 
\begin{equation}\label{eqn4.6}
\zeta_K(2n) = \frac{\tau(\chi_{D_K})(2\pi)^{4n}B_{2n}B_{2n,\bar{\chi}_{D_K}}}{4((2n)!)^2D_K^{2n}},
\end{equation}
when $K$ is real quadratic number field.

We remark that $L(s,\chi_{D_K})$ at odd positive integers is still not known. Thus the expression for Dedekind zeta function of real quadratic number field at odd positive integers can not be obtained from \eqref{eqn6.6}. Therefore we are interested to get an expression for this case too. We use Theorem \ref{mainthm} to get the same.

Let $\sigma^{(1)}(m) =  \sigma'_{-2n-2, v_K}(m)$ and $C = \frac{h_K R_K}{w_K D_K^{1/2}}$. We define
$$ P(x): = G^{0, 3}_{3, 1}\bigg(\begin{matrix}
1, x,  x \\
0, -
\end{matrix}\ \bigg|\ \frac{D_K}{4\pi ^2m}e^{-i\pi}\bigg),$$
$$Q(x):= G^{0, 3}_{3, 1}\bigg(\begin{matrix}
1, x,  x \\
0, -
\end{matrix}\ \bigg|\ \frac{D_K}{4\pi ^2m}e^{i\pi}\bigg),$$
and $$R(x):= G^{0, 3}_{3, 1}\bigg(\begin{matrix}
1, x,  x \\
0, -
\end{matrix}\ \bigg|\ \frac{D_K}{4\pi ^2m}\bigg).$$
Substite $u=2$ and $v= 2n+1$ in \eqref{realquadratic} and then apply \eqref{eqn4.6} to obtain the following expression for Dedekind zeta values of real quadratic number fields at odd integers. More precisely we have: 
\begin{theorem} For a real quadratic number field $K$ with discriminant $D_K$, we have  
\begin{eqnarray*}
&&\z_K(2n+1)=\frac{C}{\z_K(2)}\bigg(1+\frac{1}{2n}\bigg)\frac{\tau(\chi_{D_K})(2\pi)^{4n+4}B_{2n+2}B_{2n+2,\overline{\chi_{D_K}}}}{((2n+2)!)^2D_K^{2n+2}}\nonumber\\
&+&\frac{2(2\pi)^{2}}{\z_K(2)} D_K ^{- \frac{3}{2}}\sum_{m=1}^\infty m \sigma^{(1)}(m)\bigg\{ P(2)+ Q(2)+ 2 R(2)\bigg\} +\frac{2(2\pi)^{4n}}{\z_K(2)}  D_K ^{\frac{-1-4n}{2}}\nonumber\\
&\times& \sum_{m=1}^\infty m^{2n}\sigma^{(1)}(m)\bigg\{ P(2n+1)+ Q(2n+1)+ 2 R(2n+1)\bigg\}.
\end{eqnarray*}
\end{theorem} 

\subsection{Imaginary quadratic number field case}\hspace*{10cm}$ $
In this section, we discuss an expression for Dedekind zeta values attached to imaginary quadratic number fields. This expression for odd positive values is derived from the following two important results.
\begin{prop}[Ramanujan's Identity]\label{Prop4.1}
For $\alpha \beta = \pi^2$, we have
\begin{align*}
&\alpha^{-n}\left\lbrace\frac{1}{2}\zeta(2n+1) + S_{2n+1}(2\alpha) \right\rbrace = (-\beta)^{-n}\left\lbrace \frac{1}{2}\zeta(2n+1) + S_{2n+1}(2\beta) \right\rbrace \nonumber \\
& - 2^{2n}\sum_{k=0}^{n+1}(-1)^k \frac{B_{2k}B_{2n+2-2k}}{(2k)!(2n+2-2k)!}\alpha^{n+1-k}\beta^k,
\end{align*}
where $S_n(r)=\sum\limits_{k=1}^\infty \frac{1}{k^n(e^{\pi rk}-1)}$.
\end{prop} 
In \cite{CL}, M. Chamberland and P. Lopatto obtained an explicit formula for Riemann zeta odd positive values from Proposition \ref{Prop4.1}. More precisely, they obtained the following: 
\begin{prop}[Chamberland--Lopatto]
For $n= 4m-1$, 
\begin{align}\label{eq6.10}
\zeta(4m-1)&= -\frac{F_{2m-1}4^{4m-2}}{D_m}S_{4m-1}(\pi)+\frac{G_{2m-1} 4^{2m-1}}{D_m}S_{4m-1}(2\pi)\nonumber\\
&-\frac{F_{2m-1}4^{2m-1}}{D_m}S_{4m-1}(4\pi), 
\end{align}
and for $n=4m+1$, 
\begin{align}\label{eq6.11}
\zeta(4m+1) = &-\frac{16^m(G_{2m}16^m)+2E_m}{(-1+16^m)E_m}S_{4m+1}(\pi)\nonumber\\
&-\frac{2G_{2m}K_m 16^m(2 \cdot 16^m - (-4)^m +1)}{(-1+16^m)E_m}S_{4m+1}(2\pi)\nonumber\\
&-\frac{G_{2m}16^m+ 4G_{2m}16^m K_m +2E_m}{(-1+16^m)E_m}S_{4m+1}(4\pi),
\end{align}
where $$F_m=\sum\limits_{k=0}^{m+1}(-1)^k\frac{B_{2k}B_{2m+2-2k}}{(2k)!(2m+2-2k)!},$$
$$G_m=\sum\limits_{k=0}^{m+1}(-4)^k\frac{B_{2k}B_{2m+2-2k}}{(2k)!(2m+2-2k)!},$$
$$D_m=4^{2m-1}[(4^{2m-1}+1)F_{2m-1}-G_{2m-1}]/2,$$
$$K_m=\frac{\frac{1}{2}(1-4^{2m})}{1+(-4)^m-2^{4m+1}},$$
$$E_m=\frac{4^{2m}}{2}G_{2m}-2^{4m+1}K_mH_m-2^{4m}K_mG_{2m},$$
$$H_m=\sum\limits_{k=0}^m(-4)^{m+k}\frac{B_{4m}B_{4m+2-4k}}{(4k)!(4m+2-4k)!}.$$
\end{prop}

Let $\chi$ be any Dirichlet character with modulus $f$ such that $\chi(-1)=-1$. In that case H. -W. Leopoldt \cite{Leopoldt} obtained the odd $L$-values,
\begin{equation}\label{eq6.12}
L(2n+1, \chi) = (-1)^{n+1}\frac{\tau(\chi)}{2i} \left(\frac{2\pi}{f}\right)^{2n+1} \frac{B_{2n+1 ,\bar{\chi}}}{(2n+1)!},
\end{equation}
where $\tau(\chi)$ is the Gauss sum associated to $\chi$ and $B_{2n+1, \chi}$ is $(2n+1)$-th generalized Bernoulli number associated to $\chi$.

One can find the formulae for $\zeta_K(4n-1)$ and $\zeta_K(4n+1)$ by applying \eqref{eq6.10}, \eqref{eq6.11} and \eqref{eq6.12} in \eqref{eqn6.6}. 

On the other hand, $L(s,\chi_{d_K})$ at such integers is still not known. Thus an expression for Dedekind zeta values associated to imaginary quadratic number field at even positive integers can not be obtained from \eqref{eqn6.6}. Therefore it is interesting to find an expression for this case too. We use Theorem \ref{mainthm} to find such an expression. 

We put $u=2$ and $v= 2n$ in \eqref{imaginaryquadratic} and apply the above resulting expressions for $\zeta_K(4n-1)$ and $\zeta_K(4n+1)$ to get the expression for $\zeta_K(2n)$.

\begin{remark}
One can get an expression for $\zeta_K(4n-1)$ by substituting $u=v=2m$ in \eqref{imaginaryquadratic}. This expression is simpler than that of the above one, and it does not require Zagier's result for $\zeta_K(2)$. One can find a relationship between $(\zeta_K(2m+1))^2$ and $\zeta_K(4m+1)$ by substituting $u=v=2m+1$ \eqref{imaginaryquadratic}. Some other formulas can also be deduced from Theorem \ref{mainthm}. 
\end{remark}

\section{Concluding remarks}
D. Zagier in \cite{ZA1986} asked the following questions related to Theorem \ref{Zagier4}:
\begin{itemize}
\item[(Q1)] Whether there is a formula of the type Eq. \ref{zetac2} for $\z_K(4),\ \z_K(6)$, etc attached an arbitrary imaginary quadratic number field $K$?
\item[(Q2)] Whether Theorem \ref{Zagier1} can be sharpened similar to Theorem \ref{Zagier4} for any arbitrary number field? 
\end{itemize}
He also remarked that the answer to this question is possible if one can prove Theorem \ref{Zagier1} and Theorem \ref{Zagier4} using methods available in analytic number theory. Here, we have obtained the expressions for $\z_K(m)$ attached to an arbitrary quadratic number field $K$ at any positive integer ($m \geq 2$). 
We believe that these expressions can be useful to provide satisfactory answer to (Q1).
It is possible to obtain analogous expression for higher degree number fields which may be helpful to answer (Q2). 

We conclude by making the following comments: 
\begin{itemize}
\item[(1)] It could be possible to determine special values of Dedekind zeta functions for any quadratic number fields from the explicit formulae. 
But, this would involve the evaluation of Meijer G-function and other complicated terms. 
\item[(2)] One can obtain approximate functional equation for the product of two Dedekind zeta functions in the critical region by following the arguments of J. R. Wilton \cite{WI30}. 
\end{itemize}

\section*{Acknowledgments}
S. Banerjee is partially supported by Infosys grant and A. Hoque is supported by the SERB-NPDF (PDF/2017/001958), Govt. of India. 
The authors are grateful to Professors Masumi Nakajima and Bruce C. Berndt for their careful reading, helpful comments and suggestions. The authors are thankful to the anonymous referee(s) for his/her/their valuable comments and suggestions which have helped improving the presentation immensely and drawing some useful remarks.


\begin{thebibliography}{25}
\bibitem{BM14} D. Banerjee and J. Mehta, {\it Linearlized product of two Riemann zeta functions}, Proc. Japan Acad. Ser. A Math. Sci. {\bf 90} (2014), no. 8, 123--126.

\bibitem{Erdelyi} H. Bateman, {\it Higher Transcendental Functions}, Vol.1, McGraw-Hill Book Company , New York, (1953).

\bibitem{CS16} K. Chakraborty, S. Kanemitsu and H. Tsukada, {\it Vistas of special functions II}, World Scientific, Singapore, 2007. 

\bibitem{CL} M. Chamberland and P. Lopatto {\it Formulas for odd zeta values and powers of $\pi$}, J. integer seq., {\bf 14}, (2011), 1--5.

\bibitem{CM52} K. Chandrasekharan and S. Minakshisundaram, {\it Typical means}, Oxford Univ. Press, Oxford, 1952.

\bibitem{HL1921} G. H. Hardy, J. E. Littlewood, {\it The zeros of Riemann's zeta-function on the critical line}, Math. Z. {\bf 10} (1921), no. 3-4, 283--–317.

\bibitem{HR64} G. H. Hardy and M. Riesz, {\it The general theory of Dirichlet's series}, Cambridge Tracts in Mathematics and Mathematical Physics, No. 18, Stechert-Hafner, Inc., New Yorks, 1964.

\bibitem{HE1921} E. Hecke, {\it Bestimmung der Klassenzahl einer neuen Reihe von algebraischen Zahlk\"orpern}, Nachr. Akad. Wiss. G\"ottingen Math. -Phys. Kl. II (1921), 1-23.




\bibitem{LA1968} H. Lange, {\it \"Uber eine Gattung elemetar-arithmetischer Klasseninvarianten
reell-quadratischer Zahlk\"orper}, J. Reine Angew. Math. {\bf 233} (1968), 123--175.

\bibitem{Leopoldt} H.-W. Leopoldt, {\it Eine Verallgemeinerung der Bernoullischen Zahlen}, Abh. Math. Sem. Univ. Hamburg. {\bf 22} (1958), 131--140.

\bibitem{ME1957} C. Meyer, {\it \"Uber einige Anwendungen Dedekindscher Summen}, J. Reine Angew.
Math. {\bf 198} (1957), 143--203.

\bibitem{ME1966} C. Meyer, {\it \"Uber die Bildung von elementar-arithmetischen Klasseninvarianten in reell-quadratischen Zahlk\"orpern}, ``Algebraische Zahlentheorie", Mathematisches Forschungsintitut Oberwolfach, Berichte, Heft {\bf 2} (1966), 165--215.

\bibitem{NA03} M. Nakajima, {\it A new expression for the product of the two Dirichlet series, I}, Proc. Japan Acad. Ser. A Math. Sci. {\bf 79} (2003), 19--22.

\bibitem{N62} R. Narain, {\it The G-functions as unsymmetrical Fourier kernels-I}, Proc. Amer. Math. Soc. {\bf 13} (1962), 950-959.

\bibitem{SH1976} T. Shintani, {\it On the evaluation of zeta functions of totally real algebraic number fields at non-positive integers}, J. Fac. Sci. Univ. Tokyo Sect. IA Math. {\bf 23} (1976), 393--417.

\bibitem{SI1937}  C. L. Siegel, {\it \"{U}ber die analytische Theorie der quadratischen Formen. III}, Ann. of Math.  (2) {\bf 38} (1937), no. 1, 212--291.

\bibitem{SI1969} C. L. Siegel, {\it Berechung von Zetafunktionen an Ganzzahligen Stellen}, Nachr. Akad. Wiss. G\"{o}tt. \textbf{10} (1969), 87--102  (Ges. Abh. \textbf{4}, 82--97).

\bibitem{WB2017} N. L. Wang and S. Banerjee, {\it On the product of Hurwitz zeta-functions}, Proc. Japan Acad. Ser. A Math. Sci. {\bf 93} (2017), no. 5, 31--36.

\bibitem{WI30} J. R. Wilton, {\it An approximate functional equation for the product of two $\z$-functions}, Proc. London Math. Soc., {\bf S2-31} (1930), 11--17. 

\bibitem{ZA1986} D. Zagier, {\it Hyperbolic manifolds and special values of Dedekind zeta-functions}, Invent. Math. {\bf 83} (1986), 285--301.

\bibitem{ZA1976} D. Zagier, {\it On the values at negative integers of the zeta function of a real quadratic field}, Enseignement Math. {\bf 22} (1976), 55--95.

\bibitem{ZA1975} D. Zagier, {\it A Kronecker limit formula for real quadratic fields}, Math. Ann. {\bf 213} (1975), 253--284.
\end{thebibliography}
\end{document}